\documentclass{article}

\usepackage{amssymb}
\usepackage{amsmath}

\usepackage[dvips]{graphicx}
\begin{document}

%%% remove comment delimiter ('%') and select language if required
%\selectlanguage{spanish} 
\begin{center}
\noindent \textbf{                                                                                                                                                                                    A real p-homogeneous seminorm with square property is submultiplicative}
\end{center}

\noindent \textbf{}

\begin{center}
\noindent\textbf{M. El Azhari}
\end{center}
\begin{center}
Department of Mathematics, Ecole Normale Sup\'{e}rieure- Rabat, Morocco
\end{center}
\begin{center}
Email: mohammed.elazhari@yahoo.fr
\end{center}

\noindent \textbf{ } 

\noindent \textbf{Abstract.}  We give a functional representation theorem for a class of real p-Banach algebras. This theorem is used to show that every p-homogeneous seminorm with square property on a real associative algebra is submultiplicative.

\noindent \textbf{}
 

\noindent \textbf{}

\noindent \textbf{}

\noindent \textbf{1. Introduction}

\noindent \textbf{}

\noindent J. Arhippainen [1] has obtained the following result:

\noindent \textbf{}

\noindent Theorem 1 of [1]. Let q be a p-homogeneous seminorm with square property on a complex associative algebra A. Then

\begin{enumerate}
\item  Ker(q) is an ideal of A ;
\item  the quotient algebra A/Ker(q) is commutative ;
\item  q is submultiplicative ;
\item  $q^{\frac{1}{p}}$ is a submultiplicative seminorm on A.
\end{enumerate}
 
\noindent  This result is a positive answer to a problem posed in [2] and considered in [3], [4] and [5]. The proofs of (3) and (4) depend on (2) which is obtained by using a locally bounded version of the Hirschfeld-Zelazko Theorem [1, Lemma 1]. This method can not be used in a real algebra; if q is the usual norm defined on the real algebra H of quaternions, Ker (q) = \{0\} and
H/Ker (q) $\cong $ H is noncommutative, then the assertion (2) does not hold in the real case.  

\noindent  The purpose of this paper is to provide a real algebra analogue of the above Arhippainen Theorem, this improves the result in [6]. Our method is based on a functional representation theorem which we will establish; it is an extension of the Abel-Jarosz Theorem [7, Theorem 1] to real p-Banach algebras. We also give a functional representation theorem for a class of complex p-Banach algebras. As a consequence, we obtain the main result in [8].

\noindent \textbf{} 
 
\noindent \textbf{2. Preliminaries }

\noindent \textbf{  }

\noindent \textbf{  }Let A be an associative algebra over the field K=R or C. Let $ p \in ]0, 1]$, a map $\Vert.\Vert:A\to[0,\infty[$ is a p-homogeneous seminorm if for $a, b$ in A and $\alpha $ in K,  $\Vert a +b\Vert\le \Vert a \Vert+\Vert b \Vert$ and  $\Vert\alpha a\Vert=\vert\alpha\vert^{p}\Vert a\Vert$. Moreover, if  $\Vert a\Vert=0$ imply that $a = 0$, $\Vert.\Vert$ is called a p-homogeneous norm. A 1-homogeneous seminorm (resp.norm) is called a seminorm (resp.norm). $\Vert.\Vert$ is submultiplicative if  $\Vert ab\Vert\le\Vert a\Vert\Vert b\Vert$  for all $a, b$ in A. $\Vert.\Vert$ has the square property if   $\Vert a^{2}\Vert=\Vert a\Vert^{2}$ for all $a\in A$. If $\Vert.\Vert$ is a submultiplicative p-homogeneous norm on A, then $(A,\Vert.\Vert)$ is called a p-normed algebra, we denote by M (A) the set of all nonzero continuous multiplicative linear functionals on A.  A complete p-normed algebra is called a p-Banach algebra.  A uniform p-normed algebra is a p-normed algebra $(A,\Vert.\Vert)$ such that $\Vert a^{2}\Vert=\Vert a\Vert^{2}$ for all $a\in A$. Let A be a complex algebra with unit e, the spectrum of an element $a\in A$ is defined by  $Sp (a) = \{ \alpha\in C,\alpha e - a\notin A^{-1} \}$ where $A^{-1}$ is the set of all invertible elements of A.  Let A be a real algebra with unit e, the spectrum of $a\in A$ is defined by  $Sp (a) = \{s + it\in C,(a - se)^{2} + t^{2}e\notin A^{-1}  \}$.  Let A be an algebra, the spectral radius of an element $a\in A$ is defined by $r(a)=\sup\lbrace\vert\alpha\vert,\alpha\in Sp(a)\rbrace$.  Let $(A,\Vert.\Vert)$ be a p-normed algebra, the limit  $\lim_{n\to\infty }\Vert a^{n}\Vert^{\frac{1}{pn}} $  exists for each $a\in A$, and if A is complete, we have $r(a)= \lim_{n\to \infty}\Vert a^{n}\Vert^{\frac{1}{pn}} $ for all $a\in A$.  A $\ast$-algebra is a complex algebra with a mapping $\ast: A\to A, a\to a^{\ast}$ , such that, for $a, b$ in A and $\alpha\in C,(a^{\ast})^{\ast}=a,(a+b)^{\ast}=a^{\ast}+b^{\ast},(\alpha a)^{\ast}=\bar{\alpha}a^{\ast},(ab)^{\ast}=b^{\ast}a^{\ast}$. The map $\ast$ is called an involution on A. An element  a $\in$ A  is said to be hermitian if $a^{\ast}=a$ . The set of all hermitian elements of A is denoted by H (A).
 
\noindent \textbf{}
 
\noindent \textbf{3. A functional representation theorem for a class of real p-Banach algebras}

\noindent \textbf{  }
 
\noindent \textbf{  }We will need the following result due to B. Aupetit and J. Zemanek ([9] or [10]), their algebraic approach works for real p-Banach algebras.

\noindent \textbf{}

\noindent \textbf{Theorem 3.1.} Let $A$ be a real p-Banach algebra with unit. If there is a positive constant $\alpha$  such that $r(ab)\leq\alpha r(a)r(b)$ for all $a,b$ in $A$ , then for every irreducible representation $\pi$ of $A$  on a real linear space $E$, the algebra  $\pi(A)$ is isomorphic (algebraically) to its commutant in the algebra $L(E)$ of all linear transformations on $E$ .

\noindent \textbf{}

\noindent Let $A$ be a real p-Banach algebra with unit such that  $\Vert a\Vert^{\frac{1}{p}} \le m r(a)$  for some positive constant $m$ and all $a\in A$.  Let $X(A)$ be the set of all nonzero multiplicative linear functionals from A into the noncommutative algebra $H$ of quaternions. For $a\in A$, we consider the map $J(a):X(A)\to H, J(a)x = x(a)$ for all $x\in X(A)$. We endow X (A) with the weakest topology such that all the functions $J(a), a\in A$, are continuous. The map $J:A\to C(X(A),H),a\to J(a)$, is a homomorphism from $A$ into the real algebra of all continuous functions from $X(A)$ into $H$.

\noindent \textbf{}

\noindent \textbf{Theorem 3.2.} If $\pi$ is an irreducible representation of $A$, then $\pi(A)$ is isomorphic to $R$, $C$ or $H$.

\noindent \textbf{}

\noindent Proof.  Let $a,b\in A$ and $n\ge 1$, we have $\Vert (ab)^{n}\Vert\leq\Vert a\Vert^{n}\Vert b\Vert^{n} $, then $\Vert (ab)^{n}\Vert^{\frac{1}{pn}}\leq\Vert a\Vert^{\frac{1}{p}}\Vert b\Vert^{\frac{1}{p}}$ . Letting $n\to\infty $, we obtain $r(ab)\leq m^{2}r(a)r(b)$. Let $\pi $ be an irreducible representation of $A$ on a real linear space $E$. By Theorem 3.1, $\pi(A)$ is isomorphic to its commutant Q in the algebra $L(E)$ of all linear transformations on $E$. Let $y_{0} $ be a fixed nonzero element in $E$. For $y\in E$, we consider $\Vert y\Vert_{E}=\inf\lbrace\Vert a\Vert, a\in A $ and $ \pi(a)y_{0}=y\rbrace $   . By the same proof as in [11, Lemma 6.5], $\Vert.\Vert_{E}$ is a p-norm on $E$ and Q is a real division p-normed algebra of continuous linear operators on $E$. By [12], Q is isomorphic to $R$, $C\ $or $H$.

\noindent \textbf{}

\noindent \textbf{Proposition 3.3.} $A$ is semisimple and $X(A)$ is a nonempty set which separates the elements of A.

\noindent \textbf{}

\noindent Proof.  By the condition $\Vert a\Vert^{\frac{1}{p}}\leq mr(a)$ for all $a\in A$, we deduce that $A$ is semisimple. Let $a$ be a nonzero element in $A$, since $A$ is semisimple, there is an irreducible representation $\pi$ of $A$ such that $\pi(a)\ne 0$. By Theorem 3.2, there is $\varphi:\pi(A)\to H$ an isomorphism (into). We consider the map $T = \varphi o\pi, T:A\to H$ is a multiplicative linear functional. Moreover, $T(a)=\varphi(\pi(a))\ne 0 $ since $\pi(a)\ne 0$ and $\varphi $ is injective.

\noindent \textbf{}

\noindent \textbf{Proposition 3.4.}  
\begin{enumerate}
\item  $\vert x(a)\vert\leq\Vert a\Vert^{\frac{1}{p}}$ for all $a\in A$ and $x\in X(A);$ 
\item  An element $a$ is invertible in A if and only if $J(a)$ is invertible in $C(X(A),H);$
\item  $Sp(a) = Sp(J(a))$ for all $a\in A.$
\end{enumerate}

\noindent \textbf{}

\noindent Proof. (1): Since $H$ is a real uniform Banach algebra under the usual norm $\vert.\vert,    \vert x(a)\vert=r_{H}(x(a))\le r_{A}(a)\le\Vert a\Vert^{\frac{1}{p}} $   for all $a\in A$ and $x\in X(A).$

\noindent (2): The direct implication is obvious. Conversely, let $\pi $ be an irreducible representation of $A$. By Theorem 3.2, there is $\varphi:\pi(A)\to H$ an isomorphism (into). Since $\varphi o\pi\in   X(A)$ and $J(a)$ is invertible, $0 \ne J(a)(\varphi o\pi)=\varphi (\pi(a))$, then $\pi(a)\ne 0.$  Consequently, a is invertible.

\noindent  (3):   $s+it\in Sp(a)$ iff $ (a-se)^{2}+t^{2}e \notin A^{-1} $ 

\noindent Iff  $J((a-se)^{2}+ t^{2}e) \notin  C(X(A),H)^{-1}$  by (2) 

\noindent Iff  $(J(a)-sJ(e))^{2}+t^{2}J(e)\notin C(X(A),H)^{-1}$

\noindent Iff  $s+it\in Sp(J(a)).$

\noindent \textbf{}

\noindent \textbf{Proposition 3.5.} $X(A)$ is a Hausdorff compact space.

\noindent \textbf{}

\noindent Proof.  Let $ x_{1}, x_{2} $ in $ X(A), x_{1}\ne x_{2} $ , there is an element $a\in A$ such that $x_{1}(a)\neq x_{2}(a)$ , i.e.  $J(a)x_{1}\neq J(a)x_{2} $ , so  $X(A)$ is Hausdorff. Let $a\in A $ and $K_{a}= \lbrace q\in H, \vert q\vert\le\Vert a\Vert^{\frac{1}{p}}\rbrace $ , $K_{a}$ is compact in  $H$. Let  $K$ be the topological product of $K_{a}$ for all $a\in A, K$ is compact by the Tychonoff Theorem. By Proposition 3.4(1),  $X(A)$ is a subset of $K$. It is easy to see that the topology of $X(A)$ is the relative topology from  $K$ and that $X(A)$ is closed in $K$. Then $X(A)$ is compact.

\noindent \textbf{}

\noindent \textbf{Theorem 3.6.} The map $J:A\to C(X(A),H), a\to J(a),$ is an isomorphism (into) such that $m^{-1}\Vert a\Vert^{\frac{1}{p}}\le\Vert J(a)\Vert_{s}\le\Vert a\Vert^{\frac{1}{p}}$   for all $a\in A,$ where $\Vert.\Vert_{s}$  is the supnorm on $C(X(A),H).$ If $m=1$ , we have $\Vert a\Vert^{\frac{1}{p}}=\Vert J(a)\Vert_{s}$  for all $a\in A.$

\noindent \textbf{}

\noindent Proof.  By Proposition 3.3,  $J$ is an injective homomorphism. Let $a\in A,$ by Proposition 3.4(3), $r(a)=r(J(a))=\Vert J(a)\Vert_{s} $  since $C(X(A),H)$ is a real uniform Banach algebra under the supnorm $\Vert.\Vert_{s} $ . Moreover, $\Vert J(a)\Vert_{s}\leq\Vert a\Vert^{\frac{1}{p}} $  by Proposition 3.4(1). Then  $m^{-1}\Vert a\Vert^{\frac{1}{p}}\le r(a)=\Vert J(a)\Vert_{s}\le\Vert a\Vert^{\frac{1}{p}}.$  

\noindent \textbf{}

\noindent   As an application, we obtain an extension of the Kulkarni Theorem [13, Theorem 1] to real p-Banach algebras

\noindent \textbf{}

\noindent \textbf{Theorem 3.7.} Let $a$ be an element in $A$ such that $Sp(a)\subset R,$ then $a$ belongs to the center of $A.$

\noindent \textbf{}

\noindent Proof. By Theorem 3.6, $J:A\to C(X(A),H)$  is an isomorphism (into). Let $a\in A$  with $Sp(a) \subset R.$ Let $x\in X(A)$ and $x(a)=s+t$   where $s\in R$  and $t= t_{1}i+t_{2}j+t_{3}k.$  Suppose that $t\ne 0.$  We have $ (x(a)-s)^{2} =t^{2} = -(t_{1}^{2}+t_{2}^{2}+t_{3}^{2})= -\vert t\vert^{2} $, then $ (x(a)-s)^{2}+\vert t\vert^{2}= 0$. Consequently $s+i\vert t\vert \in Sp(x(a))\subset Sp(a)$  with $ \vert t\vert\ne 0,$  a contradiction. Then $J(a)\in C (X(A),R)$  and 

\noindent $ J(a)J(b)= J(b)J(a)$  for all $b$ in $A,$ i.e.  $J(ab-ba)= 0 $ for all $b$ in $A$. Since $J$ is injective, $ab-ba=0$ for all $b$ in $A.$
 
\noindent \textbf{}
 
\noindent \textbf{4. A functional representation theorem for a class of complex p-Banach algebras}

\noindent \textbf{  }
 
\noindent \textbf{  }Let $\Vert.\Vert$ be a submultiplicative p-homogeneous seminorm on a complex algebra  $A$. For $a\in A, \vert a\vert $ is defined as follows: $\vert a\vert =\inf\sum^n_{i=1}\Vert a_{i}\Vert^{\frac{1}{p}} $, where the infimum is taken over all decompositions of $a$ satisfying the condition $a=\sum^n_{i=1}a_{i}$, $ a_{1},\dots,a_{n}\in A.$ By [14, Theorem 1], $\vert.\vert$ is a submultiplicative seminorm on $A,$ it is called the support seminorm of $\Vert.\Vert $. Also, it is shown [14] the following result: 

\noindent \textbf{}

\noindent Theorem 2 of [14].  Let $A$ be a complex algebra, $ \Vert.\Vert $ a submultiplicative p-homogeneous seminorm on $A$, and $\vert.\vert $ the support seminorm of $\Vert.\Vert$ . Then $ \lim  _{n\to \infty}\Vert a^{n}\Vert^{\frac{1}{pn}}=\lim_{n\to \infty}\vert a^{n}\vert^{\frac{1}{n}}$   for all $a\in A.$

\noindent \textbf{}

\noindent In the proof of this theorem, Xia Dao-Xing uses the following inequality: If $a=a_{1}+\cdots +a_{m}$  and $n\ge 1$, then  $\Vert a^{n}\Vert\le \sum_{{\alpha }_{1}+\cdots +{\alpha }_{m} =n}      (\frac{n!}{\alpha_{1}!\cdots\alpha_{m}!})^{p}\Vert a_{1}\Vert^{\alpha_{1}}\cdots\Vert a_{m}\Vert^{\alpha_{m}} $. If the algebra is commutative, $a^{n}= (a_
{1}+\cdots +a_{m})^{n}=\sum_{\alpha_{1} +\cdots + \alpha_{m}=n}\frac{n!}{\alpha_{1}!\cdots\alpha_{m}!}     a_{1}^{\alpha_{1}}\cdots a_{m}^{\alpha_{m}}$ , then $\Vert a^{n}\Vert\le \sum_{\alpha _{1}+\cdots +\alpha_{m} =n} (\frac{n!}{\alpha_{1}!\cdots\alpha_{m}!})^{p}\Vert a_{1}\Vert^{\alpha_{1}}\cdots\Vert a_{m}\Vert^{\alpha_{m}}.$ This inequality is not justified in the noncommutative case; if the algebra is noncommutative, we only have $\Vert a^{n}\Vert\le \sum_{\alpha  _{1}+\cdots +\alpha_{m} =n}\frac{n!}{\alpha_{1}!\cdots\alpha_{m}!}\Vert a_{1}\Vert^{\alpha_{1}}\cdots\Vert a_{m}\Vert^{\alpha_{m}}$. For the sequel, we will use Theorem 2 of [14] in the commutative case.

\noindent \textbf{}

\noindent \textbf{Theorem 4.1.} Let $(A,\Vert.\Vert)$ be a complex p-normed algebra such that  $\Vert a\Vert^{2}\leq m\Vert a^{2}\Vert $ for some positive constant $m$  and all $a\in A.$ Then $\vert a\vert\leq\Vert a\Vert^{\frac{1}{p}}\leq m^{\frac{1}{p}}\vert a\vert $ and $\vert a\vert^{2}\leq m^{\frac{2}{p}}\vert a^{2}\vert$ for all $a\in A$, where  $\vert.\vert  $ is the support seminorm of $ \Vert.\Vert $ .

\noindent \textbf{}

\noindent Proof.  The completion $B$ of $(A,\Vert.\Vert )$ is a p-Banach algebra such that  $\Vert  b\Vert^{2}\leq m\Vert b^{2}\Vert $ for all $b\in B$, it is commutative by [1, Lemma 1], so $A$ is commutative. By induction, $\Vert a\Vert\le m^{1-2^{-n}}\Vert a^{2^{n}}\Vert^{2^{-n}}$  for all $a\in A$ and $n\ge 1$, then $\Vert a\Vert \le m \lim_{n\to\infty}\Vert a^{n}\Vert^{\frac{1}{n}} $  for all $a\in A.$ By the commutative version of [14, Theorem 2], we have $\vert a\vert\le\Vert a\Vert^{\frac{1}{p}}\le m^ {\frac{1}{p}}\lim_{n\to\infty}\Vert a^{n}\Vert^{\frac{1}{pn}}= m^{\frac{1}{p}}\lim_{n\to\infty}  \vert a^{n}\vert^{\frac{1}{n}}\le m^{\frac{1}{p}}\vert a\vert $ for all $a\in A.$  From the above inequalities, $\vert a\vert^{2}\le\Vert a\Vert^{\frac{2}{p}}\le (m\Vert a^{2}\Vert)^{\frac{1}{p}}\le m^{\frac{2}{p}}\vert a^{2}\vert  $.

\noindent \textbf{}

\noindent \textbf{Corollary 4.2.} Let $(A,\Vert.\Vert )$ be a complex uniform p-normed algebra. Then $ \vert a\vert =\Vert a\Vert^{\frac{1}{p}}$  for all $a\in A.$

\noindent \textbf{}

\noindent \textbf{Theorem 4.3.} Let $(A,\Vert.\Vert )$ be a complex p-Banach algebra with unit such that $\Vert a\Vert^{2}\leq m\Vert a^{2}\Vert $  for some positive constant $m$  and all $a\in A.$ Then the Gelfand map $G:A\to C(M(A))$ is an isomorphism (into) such that $m^{-\frac{2}{p}}\Vert a\Vert^{\frac{1}{p}}\le m^{-\frac{1}{p}}\vert a\vert \le\Vert G(a)\Vert_{s} \le\vert a\vert \le\Vert a\Vert^{\frac{1}{p}} $  for all $a\in A,$ where $\Vert.\Vert_{s}$ is the supnorm on $C(M(A)).$

\noindent \textbf{}

\noindent Proof.   A is commutative by [1, Lemma 1]. By Theorem 4.1, $\vert a\vert\le\Vert a\Vert^ {\frac{1}{p}}\le m^{\frac{1}{p}}\vert a\vert $ for all $a\in A$, then $(A,\vert.\vert)$ is a complex commutative Banach algebra with unit. Clearly $ M(A)= M(A,\Vert.\Vert)= M(A,\vert.\vert )$ is a nonempty compact space. As in the proof of Theorem 4.1, we have $\vert a\vert\le m^{\frac{1}{p}}\lim_{n\to\infty }  \vert a^{n}\vert^{\frac{1}{n}}= m^{\frac{1}{p}}\sup\lbrace\vert f(a)\vert, f\in M(A)\rbrace =m^{\frac{1}{p}}\Vert G(a)\Vert_{s}\le m^{\frac{1}{p}}\vert a\vert.$ Let $a\in A,$ from 

\noindent  

\noindent the above inequalities, $m^{-\frac{2}{p}}\Vert a\Vert^{\frac{1}{p}}\le m^{-\frac{1}{p}}\vert a\vert\le\Vert G(a)\Vert_{s}\le\vert a\vert\le\Vert a\Vert^{\frac{1}{p}}$.

\noindent \textbf{}

\noindent \textbf{Corollary 4.4.} Let $(A,\Vert.\Vert )$ be a complex uniform p-Banach algebra with unit. Then the Gelfand map $G:A\to C(M(A))$ is an isomorphism (into) such that $\vert a\vert =\Vert a\Vert^ {\frac{1}{p}} =\Vert G(a)\Vert_{s} $ for all $a\in A.$

\noindent \textbf{}

\noindent \textbf{Theorem 4.5.} Let $(A,\Vert.\Vert)$ be a complex p-normed $\ast$-algebra with unit such that

\begin{enumerate}
\item  $\Vert a\Vert^{2}\leq m\Vert a^{2}\Vert $ for some positive constant $m$ and all $a\in A;$

\item  Every element in $H(A)$ has a real spectrum in the completion $B$ of $A.$
\end{enumerate}

\noindent Then the involution $\ast$ is continuous on $A$ and the Gelfand map $G:B\to C(M(B))$ is a $\ast  $-isomorphism such that $m^{-\frac{2}{p}}\Vert b\Vert^{\frac{1}{p}}\le\Vert G(b)\Vert_{s}\le\Vert b\Vert^ {\frac{1}{p}}$  for all $b$ in $B.$

\noindent \textbf{}

\noindent Proof.  By Theorem 4.3, it remains to show that the involution $\ast$ is continuous on $A$, $G(b^{\ast})=G(b)^{\ast}$ for all $b\in B,$ and $G$ is surjective. Let $h\in H(A), Sp_{B}(h)=\lbrace f(h), f\in M(B)\rbrace   \subset R$  by (2). Let $a\in A,$ we have $a=h_{1}+ih_{2}$  with $h_{1} ,h_{2}\in H(A).$ Let $f\in M(B),
f(a^{\ast})=f(h_{1}-ih_{2})=f(h_{1})-if(h_{2})=(f(h_{1})+if(h_{2}))^{\ast}=f(h_{1}+ih_{2})^{\ast}= f(a)^{\ast}$ since $f(h_{1})$ and $f(h_{2})$ are real.  
Then $G(a^{\ast})=G(a)^{\ast} $ for all $a\in A.$ By Theorem 4.3, $m^{-\frac{2}{p}}\Vert a^{\ast} \Vert^{\frac{1}{p}}\le\Vert G(a^{\ast})\Vert_{s}=\Vert G(a)^{\ast}\Vert_{s}=\Vert G(a)\Vert_{s}\le\Vert a\Vert^{\frac{1}{p}}$  for all $a\in A,$ then $ \Vert a^{\ast}\Vert\le m^{2}\Vert a\Vert $  for all $a\in A.$ Consequently, the involution $\ast$ is continuous on A and can be extended to a continuous involution on $B$ which we will also denote by $\ast.$ Let $b\in B,$ there exists a sequence $(a_{n})_{n}$ in $A$  such that $a_{n}\to b$. Since the involution on $B$ and the Gelfand map $G:B\to C(M(B))$ are continuous, we have $G(a^{\ast}_{n})\to G(b^{\ast})$ and $G(a_{n})^{\ast}\to\ G(b)^{\ast},$ then $G(b^{\ast}) =G(b)^{\ast} .$ By the Stone-Weierstrass Theorem, we deduce that $G$ is surjective.

\noindent \textbf{}

\noindent   As a consequence, we obtain the main result in [8].

\noindent \textbf{}

\noindent \textbf{Corollary 4.6.} Let $A$ be a complex uniform p-normed $\ast$-algebra with unit such that every element in $H(A)$ has a real spectrum in the completion $B$ of $A.$ then $B$ is a commutative $C^{\ast}$-algebra.

\noindent \textbf{}

\noindent \textbf{5. The main result}

\noindent \textbf{}
 
\noindent \textbf{Theorem 5.1.} Let $A$ be a real associative algebra. Every p-homogeneous seminorm $q$  with square property on $A$  is submultiplicative and $q^{\frac{1}{p}}$ is a submultiplicative seminorm on $A.$

\noindent \textbf{}

\noindent Proof.  By [1], there exists a positive constant $m$ such that $q(ab)\le mq(a)q(b)$ for all $a,b\in A$. $Ker(q)$ is an ideal of $A,$ the norm $ \vert.\vert $  on the quotient algebra $A/Ker(q)$ defined by $ \vert a+Ker(q)\vert=q(a)$  is a p-norm with square property. Define $ \Vert a+Ker(q)\Vert=m \vert a+Ker(q)\vert $ for all $a\in A.$ Let $a,b\in A, \Vert ab+Ker(q)\Vert=m\vert ab+Ker(q)\vert\le m^{2} \vert a+Ker(q)\vert\vert b+Ker(q)\vert=\Vert a+Ker(q)\Vert\Vert b+Ker(q)\Vert ,$ then $(A/Ker(q),\Vert.\Vert)$ is a real p-normed algebra. Let $a\in A, \Vert a^{2}+Ker(q)\Vert =m \vert a^{2}+Ker(q)\vert =m \vert a+Ker(q)\vert^{2}= m^{-1}(m\vert a+Ker(q)\vert)^{2}= m^{-1}\Vert a+Ker(q)\Vert^{2}$ i.e. $ \Vert a+Ker(q)\Vert^{2}= m\Vert a^{2}+Ker(q)\Vert.$ The completion $B$ of $(A/Ker(q),\Vert.\Vert )$ satisfies   
also the property $\Vert b\Vert^{2}= m\Vert b^{2}\Vert $ for all $b\in B,$ and by induction $\Vert b\Vert = m^{1-2^{-n}}\Vert b^{2^{n}}\Vert^{2^{-n}}$ for all $b\in B$ and $n\ge 1,$ then 
$\Vert b\Vert= mr(b)^{p}$ for all $b\in B.$ We consider two cases:
 
\noindent $B$ is unital:  By section 3, $X(B)$ is a nonempty compact space and the map $J:B\to C(X(B),H)$ is an isomorphism (into). By Proposition 3.4(3), $r(b)=r(J(b))$ for all $b\in B$. Let $b\in B, \Vert b\Vert=mr(b)^{p}=mr(J(b))^{p}=m\Vert J(b)\Vert_{s}^{p}$  since $C(X(B),H)$ is a real uniform Banach algebra under the supnorm $\Vert.\Vert_{s}$. Then $\vert b\vert =m^{-1}\Vert b\Vert =\Vert J(b)\Vert_{s}^{p} $ for all $b\in A/Ker(q)$, so $\vert.\vert $ is submultiplicative and $\vert.\vert^ {\frac{1}{p}} $  is a submultiplicative norm. Consequently, $q$ is submultiplicative and $q^{\frac{1}{p}}  $ is a submultiplicative seminorm.

\noindent $B$ is not unital:  Let $B_{1}$ be the algebra obtained from $B$ by adjoining the unit. By the same proof of [15, Lemma 2] which works for real p-Banach algebras, there exists a p-norm $N$ on $B_{1}$ such that

\begin{enumerate}
\item  $ (B_{1},N)$ is a real p-Banach algebra with unit;
\item  $ N(b)^{\frac{1}{p}}\le m^{3}r_{B_{1}}(b)$  for all $b\in B_{1}$;
\item  $N$ and $\Vert.\Vert $ are equivalent on $B.$
\end{enumerate}

\noindent By section 3, $X(B_{1})$ is a nonempty compact space and the map $J:B_{1}\to C(X(B_{1}),H)$ is an isomorphism (into). Let $b\in B$,

\noindent $\Vert b\Vert =mr_{B}(b)^{p}=mr_{B_{1}}(b)^{p}$    by (3)

\noindent                             $=mr(J(b))^{p}$    by Proposition 3.4(3)  

\noindent                             $=m\Vert J(b)\Vert_{s}^{p}$ by the square property of the supnorm.  

\noindent Then $\vert b\vert=m^{-1}\Vert b\Vert=\Vert J(b)\Vert_{s}^{p} $ for all $b\in A/Ker(q),$ so $\vert.\vert $  is submultiplicative and $\vert.\vert^{\frac{1}{p}}$ is a submultiplicative norm. Consequently, $q$ is submultiplicative and $q^{\frac{1}{p}}$ is a submultiplicative seminorm.

\noindent \textbf{}   

\noindent \textbf{}  
   
\noindent \textbf{References}

\noindent \textbf{}

\noindent [1] J. Arhippainen, On locally pseudoconvex square algebras, Publicacions Matematiques, 39 (1995), 89-93.
 
\noindent [2] S. J. Bhatt and D. J. Karia, Uniqueness of the uniform norm with an application to topological algebras, Proc. Amer. Math. Soc., 116 (1992), 499-504.

\noindent [3] S. J. Bhatt, A seminorm with square property on a Banach algebra is submultiplicative, Proc. Amer. Math. Soc., 117 (1993), 435-438.
  
\noindent [4] H. V. Dedania, A seminorm with square property is automatically submultiplicative, Proc. Indian. Acad. Sci. (Math. Sci), 108 (1998), 51-53.

\noindent [5] Z. Sebestyen, A seminorm with square property on a complex associative algebra is submultiplicative, Proc. Amer. Math. Soc., 130 (2001), 1993-1996.

\noindent [6] M. El Azhari, A real seminorm with square property is submultiplicative, Indian J. Pure Appl.Math. 43 (2012), 303-307.

\noindent [7] M. Abel and K. Jarosz, Noncommutative uniform algebras, Studia Math., 162 (2004), 213-218.

\noindent [8] A. El Kinani, On uniform hermitian p-normed algebras, Turk J. Math., 30 (2006), 221-231. 

\noindent [9] B. Aupetit and J. Zemanek, On the real spectral radius in real Banach algebras, Bull. Acad. Polon. Sci. S\'{e}r. Sci. Math. Astronom. Phys., 26 (1978), 969-973.

\noindent [10] J. Zemanek, Properties of the spectral radius in Banach algebras, Banach Center Publications,      8 (1982), 579-595.

\noindent [11] A. M. Sinclair, Automatic continuity of linear operators, Cambridge University Press, Cambridge 1976.

\noindent [12]  Ph. Turpin,  Sur une classe d'alg\`{e}bres topologiques, C. R. Acad. Sc. Paris, S\'{e}rie A, 263 (1966), 436-439.

\noindent [13] S. H. Kulkarni, Representations of a class of real $B^{\ast}$-algebras as algebras of quaternion-valued functions, Proc. Amer. Math. Soc., 116 (1992), 61-66.
 
\noindent [14] Xia Dao-Xing, On locally bounded topological algebras, Acta Math. Sinica, 14 (1964), 261-276.

\noindent [15] R. A. Hirschfeld and W. Zelazko, On spectral norm Banach algebras, Bull. Acad. Polon. Sci.  S\'{e}r. Sci. Math. Astronom. Phys., 16 (1968), 195-199.
 
\noindent \textbf{}
 
\end{document}